\documentclass[11pt]{amsart} 
\usepackage{amssymb,latexsym,enumerate,graphicx,bbm,mathptmx,microtype}

\newtheorem{theorem}{Theorem} 
\newtheorem{corollary}[theorem]{Corollary}

\newtheorem{exam}{Example}

\renewcommand\emptyset{\varnothing}
\def\th{$^{\text{th}}$}

\newcommand\commentout[1]{}
\newcommand\Def[1]{\emph{#1}}

\newcommand\link{\operatorname{link}}

\newcommand\height{\operatorname{height}} 
\newcommand\cone{\operatorname{cone}} 
\newcommand\ehr{\operatorname{ehr}} 
\newcommand\Ehr{\operatorname{Ehr}}

\def\Z{\mathbb{Z}}
\def\R{\mathbb{R}}

\def\B{\mathcal{B}}

\def\K{\mathcal{K}}

\def\P{\mathcal{P}}
\def\Q{\mathcal{Q}}

\def\A{\mathbf{A}}

\def\m{\mathbf{m}}
\def\v{\mathbf{v}}

\def\x{\mathbf{x}}

\def\0{\mathbf{0}}
\def\1{\mathbf{1}}

\def\h{h^*}

\begin{document}

\title{Stanley's Major Contributions to Ehrhart Theory}

\author{Matthias Beck}
\address{Department of Mathematics\\
         San Francisco State University\\
         San Francisco, CA 94132\\
         U.S.A.}
\email{mattbeck@sfsu.edu}
\urladdr{http://math.sfsu.edu/beck}

\dedicatory{Dedicated to Richard Stanley on the occasion of his 70\th birthday}

\begin{abstract}
This expository paper features a few highlights of Richard Stanley's extensive work in Ehrhart theory, the study of integer-point enumeration in rational polyhedra.
We include results from the recent literature building on Stanley's work, as well as several open problems.
\end{abstract}



\date{13 September 2015}

\thanks{We thank Ben Braun, Martin Henk, Richard Stanley, and Alan Stapledon for helpful discussions and suggestions.
This work was partially supported by the U.\ S.\ National Science Foundation (DMS-1162638).}

\maketitle


\section{Introduction}

This expository paper features a few highlights of Richard Stanley's extensive work in Ehrhart theory, with some pointers to the recent literature and open problems.
Pre-Stanley times saw two major results in this area: In 1962, Eug\`ene Ehrhart established the following fundamental theorem for a \Def{lattice polytope}, i.e., the convex hull of finitely many integer points in $\R^d$.

\begin{theorem}[Ehrhart \cite{ehrhartpolynomial}]\label{thm:ehrhart}
If $\P \subset \R^d$ is a lattice polytope and $n \in \Z_{ >0 }$ then
\[
  \ehr_\P(n) := \# \left( n \P \cap \Z^d \right) 
\]
evaluates to a polynomial in $n$ (the \Def{Ehrhart polynomial} of $\P$).
Equivalently, the accompanying generating function (the \Def{Ehrhart series} of $\P$) evaluates to a rational function:
\[
  \Ehr_\P(x) := 1 + \sum_{ n>0 } \ehr_\P(n) \, x^n = \frac{ \h_\P(x) }{ (1-x)^{ \dim(\P)+1 } } 
\]
for some polynomial $\h_\P(x)$ of degree at most $\dim(\P)$, the \Def{Ehrhart $h$-vector} of $\P$.\footnote{The Ehrhart $h$-vector is also known by the names of \emph{$\h$-vector/polynomial} and \emph{$\delta$-vector/polynomial}.}
\end{theorem}

We remark that the step from an Ehrhart polynomial to its rational generating function is a mere change of variables: the coefficients of $\h_\P(x)$ express $\ehr_\P(n)$ in the
binomial-coefficient basis $\binom n k$, $\binom {n+1} k$, \dots, $\binom {n+k} k$, where $k = \dim(\P)$.

In 1971, I.G.\ Macdonald proved the following \emph{reciprocity theorem}, which had been conjectured (and proved for several special cases) by Ehrhart.

\begin{theorem}[Ehrhart--Macdonald \cite{macdonald}]\label{thm:macdon}
The evaluation of the Ehrhart polynomial of $\P$ at a negative integer yields
\[
  \ehr_\P(-n) = (-1)^{ \dim(\P) } \ehr_{ \P^\circ } (n) \, ,
\]
where $\P^\circ$ denotes the (relative) interior of~$\P$.
Equivalently, we have the following identity of rational functions:
\[
  \Ehr_\P(\tfrac 1 x) = (-1)^{ \dim(\P) + 1 } \Ehr_{ \P^\circ } (x) \, ,
\]
where $\Ehr_{ \P^\circ } (x) := \sum_{ n>0 } \ehr_{ \P^\circ }(n) \, x^n$.
\end{theorem}

Stanley made several fundamental contributions to Ehrhart theory, starting in the 1970s. This paper attempts to highlight some of them, roughly in historical order.
The starting point, in Section \ref{sec:anand}, is Stanley's proof of the Anand--Dumir--Gupta conjecture. Section \ref{sec:reciprocity} features a reciprocity theorem
of Stanley that generalizes Theorem~\ref{thm:macdon}, and Section \ref{sec:inequ} contains Stanley's inequalities on Ehrhart $h$-vectors. Throughout the paper we mention open problems, and Section \ref{sec:beyond} gives some recent results in Ehrhart theory building on Stanley's work.


\section{``I am grateful to K.\ Baclawski for calling my attention to the work of Eug\`ene Ehrhart...''}\label{sec:anand}

The starting point for Stanley's work in Ehrhart theory is arguably his proof of the \emph{Anand--Dumir--Gupta conjecture} \cite{ananddumirgupta}; this proof appeared in a 1973 paper from which the quote of the section header is taken \cite[p.~631]{stanleymagic}. 
(One could, in fact, argue that \emph{$P$-partitions}, which were introduced in Stanley's Ph.D.\ thesis
\cite{stanleythesis} and are featured in Gessel's article in this volume, already show an
Ehrhart-theoretic flavor, but this geometric realization came slightly later in Stanley's work.)
The conjecture concerns the counting function $H_n(r)$, the number of $(n \times n)$-matrices with nonnegative integer entries that sum to $r$ in every row and
column; these matrices are often referred to as \emph{semimagic squares}. The function $H_n(r)$ goes back to MacMahon \cite{macmahon}, who computed the first
nontrivial case, $H_3 (r) = \binom {r+5} 5 - \binom {r+2} 5$. Anand, Dumir, and Gupta conjectured that
\begin{itemize}
  \item $H_n(r)$ is a polynomial in $r$ for any fixed $n$,
  \item this polynomial has roots at $-1, -2, \dots, -n+1$, and 
  \item it satisfies the symmetry relation $H_n(-r) = (-1)^{ n-1 } H_n(r-n)$. 
\end{itemize}
Stanley showed that one could use what is now called the \emph{Elliott--MacMahon algorithm} (going back to \cite{elliott} and \cite{macmahon}; see also
\cite{andrewspa1} and its dozen follow-up papers) in connection with the \emph{Hilbert syzygy theorem} (see, e.g., \cite{eisenbudbook}) to prove the Anand--Dumir--Gupta conjecture.\footnote{Stanley's comment in \cite[p.~6]{stanleyhowupperbound} is amusing in this context: ``I had taken a course in graduate school on commutative algebra that I did not find very interesting. It did not cover the Hilbert syzygy theorem. I had to learn quite a bit of commutative algebra from scratch in order to understand the work of Hilbert.''}
Stanley realized that $H_n(r)$ is closely related to the geometry of the \emph{Birkhoff--von Neumann polytope} $\B_n$, the set of all doubly-stochastic $(n \times n)$-matrices. In modern language, $H_n(r)$ equals the Ehrhart polynomial of $\B_n$, and because $\B_n$ is \emph{integrally closed} (essentially by the Birkhoff--von Neumann theorem that the extreme points of $\B_n$ are precisely the permutation matrices; see, e.g.,
\cite{brunsgubeladzektheory} for more on integrally closed polytopes), this Ehrhart series equals the \emph{Hilbert series} of the semigroup algebra generated by the integer point in $\B_n \times \{ 1 \}$, graded by the last coordinate. Stanley proved that this Hilbert series has the following properties, which are a direct translation of the Anand--Dumir--Gupta conjecture into the language of generating functions:

\begin{theorem}[Stanley \cite{stanleymagic}]\label{thm:stanleysemimagic}
For each $n$ there exists a palindromic polynomial $h_n(x)$ of degree $n^2 - 3n + 2$ such that
\[
  1 + \sum_{ r>0 } H_n(r) \, x^r = \frac{ h_n(x) }{ (1-x)^{ n^2 - 2n + 2 } } \, .
\]
\end{theorem}

One could skip the detour through commutative algebra and directly realize this rational generating
function  as the Ehrhart series of $\B_n$, with the needed properties to affirm the
Anand--Dumir--Gupta conjecture. Nevertheless, the algebraic detour is worth taking; aside from its
inherent elegance, it allowed Stanley to realize that the monomials $z_1^{ m_1 } z_2^{ m_2 } \cdots
z_{n^2}^{ m_{n^2} }$, where $(m_1, m_2, \dots, m_{ n^2 }) \in \Z_{ \ge 0 }^{ n^2 }$ ranges over all
semimagic squares, generate a \emph{Cohen--Macaulay algebra}, and this implies, following Hochster's work~\cite{hochster}:

\begin{theorem}[Stanley \cite{stanleymagiccohenmac}]\label{thm:stanleysemimagicintegral}
The polynomial $h_n(x)$ in Theorem~\ref{thm:stanleysemimagic} has nonnegative coefficients.
\end{theorem}

Stanley conjectured that the coefficients of $h_n(x)$ in Theorem~\ref{thm:stanleysemimagic} are also \emph{unimodal}
(the coefficients increase up to some point and then decrease), which was proved by Athanasiadis some three decades later \cite{athanasiadismagic}; we will say more about this in Section \ref{sec:beyond}.

Another problem connected to Theorem~\ref{thm:stanleysemimagic} (mentioned by Stanley but certainly older than his work) is still wide open, namely, the quest
for the volume of $\B_n$, which equals $h_n(1)$ after a suitable normalization. This volume is known only for $n \le 10$, though there has been recent progress, e.g.,
in terms of asymptotic and combinatorial formulas~\cite{beckpixton,canfieldmckay,deloeraliuyoshida}.

It is worth noting that Stanley proved versions of Theorems~\ref{thm:stanleysemimagic} and
\ref{thm:stanleysemimagicintegral} for the more general \emph{magic labellings} of graphs (and semimagic squares
correspond to such labellings for a complete bipartite graph $K_{ nn }$). In this more general context, the
associated counting functions become \emph{quasipolynomials} of period 2 (see, e.g., \cite[Chapter~4]{stanleyec1} for more about quasipolynomials), foreshadowing in
some sense Zaslavsky's work on enumerative properties of \emph{signed graphs}: Stanley's magic labellings are essentially flows on all-negative signed
graphs~\cite{nnz,zaslavskyorientationsignedgraphs}.

Stanley's work on the Anand--Dumir--Gupta conjecture was not just the starting point of his contributions to Ehrhart theory. As he mentions in
\cite{stanleyhowupperbound}, it opened the door to what are now called \emph{Stanley--Reisner rings} and the use of Cohen--Macaulay algebras in geometric
combinatorics, famously leading to Stanley's proof of the \emph{upper bound conjecture} for spheres \cite{stanleyupperbound}.
Stanley's appreciation for the polynomials $H_n(r)$ is also evident in his writings: they are prominently featured in his influential books \cite{stanleycombcommalg,stanleyec1}.
We close this section by mentioning that a highly readable account of commutative-algebra concepts behind the Anand--Dumir--Gupta conjecture can be found in~\cite{brunssemimagic}.


\section{Reciprocity}\label{sec:reciprocity}

Theorem~\ref{thm:macdon} is an example of a \emph{combinatorial reciprocity theorem}: we get interesting information out of a counting function when we evaluate it at a \emph{negative} integer (and so, a priori the counting function does not make sense at this number).
We remark that the last part of the Anand--Dumir--Gupta conjecture follows from Theorem~\ref{thm:macdon} applied to the Birkhoff--von Neumann polytope
$\B_n$, as one can easily show that $\ehr_{ \B_n^\circ } (r) = H_n(r-n)$.

Stanley had discovered reciprocity theorems for $P$-partitions and order polynomials in his thesis
\cite{stanleythesis}, and so it was natural for him to realize Theorem~\ref{thm:macdon} as a special case of a wider phenomenon.
A \Def{rational cone} is a set of the form $\left\{ \x \in \R^d : \, \A \, \x \le \0 \right\}$ for some integer
matrix $\A$. It is not hard to see that the multivariate generating function
\[
  \sigma_\K \left( z_1, z_2, \dots, z_d \right) := \sum_{ \left( m_1, m_2, \dots, m_d \right) \in \K \cap \Z^d }
z_1^{ m_1 } z_2^{ m_2 } \cdots z_d^{ m_d } 
\]
evaluates to a rational function in $z_1, z_2, \dots, z_d$ if $\K$ is a rational cone (see, e.g.,
\cite[Chapter~3]{ccd}).
Stanley proved that this rational function satisfies a reciprocity theorem.

\begin{theorem}[Stanley \cite{stanleyreciprocity}]\label{thm:stanleyrec}
If $\K$ is a rational cone then
\[
  \sigma_\K \left( \tfrac{ 1 }{ z_1 } , \tfrac{ 1 }{ z_2 } , \dots, \tfrac{ 1 }{ z_d } \right) = (-1)^{\dim(\K)} \, \sigma_{
\K^\circ } \left( z_1, z_2, \dots, z_d \right) .
\]
\end{theorem}

A simple proof of Theorem~\ref{thm:stanleyrec}, based on the ideas of \cite{bsstanleyirrational}, can be found in
\cite[Chapter~4]{ccd}.
True to the theme of this survey, our description of Theorem~\ref{thm:stanleyrec} is geometric, while Stanley's viewpoint presented in \cite{stanleyreciprocity} is based on integral solutions of a system of integral linear equations; in this language, the reciprocity is between \emph{nonnegative} and \emph{positive} solutions.
(It is not hard to see that both viewpoints are equivalent.)
This language also connects once more to the Elliott--MacMahon algorithm mentioned in Section~\ref{sec:anand},
and Stanley uses this connection in \cite{stanleyreciprocity} to present his \emph{monster reciprocity theorem},
which (a bit oversimplified) can be thought of an affine version of Theorem~\ref{thm:stanleyrec}. It has been
recently revitalized by Xin~\cite{xinmonsterrec}.

We finish this section with a sketch how Theorem~\ref{thm:macdon} follows as a corollary of Theorem~\ref{thm:stanleyrec}.
Given a lattice polytope $\P \subset \R^d$, we consider its \Def{homogenization}
\[
  \cone(\P) := \sum_{ \v \text{ vertex of } \P } \R_{ \ge 0 } (\v, 1)
\]
by lifting the vertices of $\P$ into $\R^{ d+1 }$ to the hyperplane $x_{ d+1 } = 1$ and taking the nonnegative
span of this lifted version of $\P$. Thus (by the \emph{Minkowski--Weyl theorem}---see, e.g., \cite[Lecture
1]{ziegler}) $\cone(\P)$ is a rational cone, and because $\cone(\P) \cap \{ \x \in \R^d : \, x_{ d+1 } = n \}$ is identical to $n \P$ (embedded in $\{ \x \in \R^d : \, x_{ d+1 } = n \}$),
\[
  \Ehr_\P(x) = \sigma_{ \cone(\P) } (1, 1, \dots, 1, x) \, .
\] 
Theorem~\ref{thm:macdon} follows now by specializing all but one variable in Theorem~\ref{thm:stanleyrec}.


\section{Ehrhart inequalities}\label{sec:inequ}

Theorem~\ref{thm:stanleysemimagicintegral} generalizes to all lattice polytopes, and this is arguably Stanley's most important contribution to the intrinsic study of Ehrhart polynomials.

\begin{theorem}[Stanley \cite{stanleydecomp}]\label{thm:nonneg}
For any lattice polytope $\P$, the Ehrhart $h$-vector $\h_\P(x)$ has nonnegative coefficients.
\end{theorem}

Even though Stanley explicitly states this result first in \cite{stanleydecomp} (and gives a proof using a shelling triangulation argument), he attributes Theorem
\ref{thm:nonneg} to \cite[Proposition 4.5]{stanleymagic}---the magic-graph-labelling version of Theorem~\ref{thm:stanleysemimagicintegral}.
Theorem~\ref{thm:nonneg} can be viewed as a starting point for the problem of classifying Ehrhart polynomials/Ehrhart $h$-vectors. This problem is wide open, already in dimension 3. (In dimension 2, it is essentially solved by Pick's theorem \cite{pick} and an inequality of Scott~\cite{scott}.)

As part of his study of monotonicity of $h$-vectors of Cohen--Macaulay complexes, Stanley deduced the following refinement of Theorem~\ref{thm:nonneg}.

\begin{theorem}[Stanley \cite{stanleymonotonicity}]\label{thm:monoton}
If $\P \subseteq \Q$ are lattice polytopes then $\h_\P(x) \le \h_\Q(x)$ (component-wise).
\end{theorem}

Theorem~\ref{thm:nonneg} can be realized as a corollary of Theorem~\ref{thm:monoton} by choosing $\P$ to be a \emph{unimodular simplex} (i.e., a $d$-dimensional lattice polytope of volume $\frac 1 {d!}$, which comes with the Ehrhart $h$-vector $\h_\P(x) = 1$).

One can show that Theorems~\ref{thm:nonneg} and \ref{thm:monoton} also hold for \emph{rational} $d$-polytopes $\P$ and $\Q$ (i.e., polytopes whose vertices have rational coordinates) if their Ehrhart series are written in the form
\[
  \frac{ \h_{\P/\Q} (x) }{ (1-x^p)^{ d+1 } } 
\]
for some $p \in \Z_{ >0 }$ such that $p\P$ and $p\Q$ are lattice polytopes.
(The accompanying Ehrhart counting functions for $\P$ and $\Q$ are then quasipolynomials, and $p$ is a period of them.)
The arguably simplest proofs of (rational versions of) Theorems~\ref{thm:nonneg} and \ref{thm:monoton} are given in~\cite{bsstanleyirrational}.

A natural question is whether there are any natural \emph{upper} bounds complementing Theorem~\ref{thm:nonneg}. Of course, the volume (and therefore $\h_\P(1)$, the sum
of the Ehrhart-$h$ coefficients) of a lattice polytope $\P$ can be arbitrarily large, but one can ask for upper bounds given certain data.

\begin{theorem}[Haase--Nill--Payne \cite{haasenillpayne}]
The volume of a lattice polytope $\P$ (and therefore also the coefficients of $\h_\P(x)$) is bounded by a
number that depends only on the degree and the leading coefficient of~$\h_\P(x)$.
\end{theorem}

This result was conjectured by Batyrev \cite{batyrevhvector} and improves a classic theorem of Lagarias and Ziegler \cite{lagariasziegler} that the volume of a
lattice polytope is bounded by a number depending only on its dimension and the number of its interior lattice points, if the latter is positive.

Theorem~\ref{thm:nonneg} can be extended in a direction different from that of Theorem~\ref{thm:monoton}, namely, one can establish inequalities \emph{among} the
coefficients of an Ehrhart $h$-vector.  Stanley derived one set of such inequalities as a corollary of his study of Hilbert functions of semistandard graded
Cohen--Macaulay domains \cite{stanleyinequ}; thus the following result holds in a more general situation. 
The \Def{degree} of a lattice polytope $\P$ is the degree of its Ehrhart $h$-vector~$\h_\P(x)$.

\begin{theorem}[Stanley \cite{stanleyinequ}]\label{thm:stanlineq}
If $\P$ is a $d$-dimensional lattice polytope of degree $s$ then its Ehrhart $h$-vector $\h_\P(x) = \h_s x^s + \h_{ s-1 } x^{ s-1 } + \dots + \h_0$ satisfies
\[
  \h_0 + \h_1 + \dots + \h_j \le \h_s + \h_{ s-1 } + \dots + \h_{ s-j } \quad \text{ for } \quad 0 \le j \le d.
\]
\end{theorem}

\noindent
(In the above theorem and below, we define $\h_j = 0$ whenever $j < 0$ or $j$ is larger than the degree of $\P$.)
Theorem~\ref{thm:stanlineq} complements inequalities discovered by Hibi around the same time \cite{hibiehrhartineq,hibilowerbound}; they were more recently improved
by Stapledon. Together with Theorem~\ref{thm:stanlineq} and the trivial inequality $\h_1 \ge \h_d$ (which follows from the facts $\h_1 = \# (\P \cap \Z^d) - d - 1$
and $\h_d = \# (\P^\circ \cap \Z^d)$), the following result gives the state of the art in terms of linear constraints for the Ehrhart coefficients that can be easily written down in general.

\begin{theorem}[Stapledon \cite{stapledondelta}]\label{thm:stapledonineq}
If $\P$ is a $d$-dimensional lattice polytope of degree $s$ and codegree $l := d+1-s$, then its Ehrhart $h$-vector $\h_\P(x) = \h_s x^s + \h_{ s-1 } x^{ s-1 } + \dots + \h_0$ satisfies
\begin{align*}
  \h_2 + \h_3 + \dots + \h_{ j+1 } \ge \h_{ d-1 } + \h_{ d-2 } + \dots + \h_{ d-j } \quad &\text{ for } \quad 0
\le j \le \lfloor \tfrac d 2 \rfloor - 1, \\
  \h_{ 2-l } + \h_{ 3-l } + \dots + \h_1 \le \h_j + \h_{ j-1 } + \dots + \h_{ j-l+1 } \quad &\text{ for } \quad 2 \le j \le d - 1. 
\end{align*}
\end{theorem}

We will say more about this theorem in Section \ref{sec:beyond} and finish this section with the remark that a unimodular $d$-simplex satisfies each of the inequalities of Theorems~\ref{thm:stanlineq} and \ref{thm:stapledonineq} with equality.


\section{Stanley \& beyond}\label{sec:beyond}

The Ehrhart $h$-vector is philosophically close to the $h$-vector of a simplicial complex. This statement can be made much more precise, as we will show in this final
section, in which we give a flavor of recent results that build on Stanley's work in Ehrhart theory. 

The starting point can once more be found in Stanley's papers; the following result follows essentially from the definition
\[
  h_T(z) := \sum_{ k=-1 }^{ d } f_k \, z^{ k+1} \, (1-z)^{ d-k }
\]
of the \Def{$h$-vector} of a given triangulation $T$ of a $d$-dimensional polytope (here $f_k$ denotes the number of $k$-simplices in $T$) and the fact that a
unimodular simplex has a trivial Ehrhart $h$-vector (a triangulation is \Def{unimodular} if all of its simplices are). Nevertheless, the following identity is an
important base case for structural properties of Ehrhart $h$-vectors to come.

\begin{theorem}[Stanley \cite{stanleydecomp}]\label{thm:unimodularehrharth}
If $\P$ is a lattice polytope that admits a unimodular triangulation $T$ then
\[
  \Ehr_\P(z) = \frac{ h_T(z) }{ (1-z)^{ \dim(\P) + 1 } } \, .
\]
In words, the Ehrhart $\h$-vector of $\P$ is given by the $h$-vector of~$T$.
\end{theorem}

The hope is now to use properties of the $h$-vector of a triangulation to say something about Ehrhart $h$-vectors; for example, if $T$ is the cone over a boundary
triangulation of $\P$ then $h_T(z)$ satisfies the \emph{Dehn--Sommerville equations} (see, e.g., \cite{deloerarambausantos} for more about triangulations).
Unfortunately, not all lattice polytopes admit unimodular triangulations in dimension $\ge 3$---in fact, most do not---and so Theorem~\ref{thm:unimodularehrharth}
needs some tweaking before we can apply it to general lattice polytopes.
This tweaking, due to Betke and McMullen \cite{betkemcmullen}, has two main ingredients: the \Def{link} of a simplex $\Delta$ in a triangulation $T$
\[
  \link(\Delta) := \left\{ \Omega \in T : \, \Omega \cap \Delta = \emptyset, \ \Omega \subseteq \Phi \text{ for some } \Phi \in T \text{ with } \Delta \subseteq \Phi
\right\} , 
\]
and its \emph{box polynomial}
\[
  B_\Delta(x) := \sum_{ \m \in \Pi(\Delta) \cap \Z^{ d+1 } } x^{ \height(\m) } 
\]
where we define
\[
  \Pi(\Delta) := \left\{ \sum_{ \v \text{ vertex of } \Delta } \!\!\!\!\lambda_\v (\v, 1) : \, 0 < \lambda_\v < 1 \right\}
\]
and $\height(\m)$ denotes the last coordinate of $\m$.
(Geometrically, $\Pi(\Delta)$ is the open fundamental parallelepiped of $\cone(\Delta)$.)
For the empty simplex $\emptyset$ of a triangulation, we set $B_\emptyset(x) = 1$.

\begin{theorem}[Betke--McMullen \cite{betkemcmullen}]\label{thm:betkemcmullen}
Fix a triangulation $T$ of the lattice polytope $\P$. Then
\[
  \h_\P(x) = \sum_{ \Delta \in T } h_{ \link(\Delta) } (x) \, B_\Delta(x) \, .
\]
\end{theorem}

If a simplex $\Delta \in T$ is unimodular then $B_\Delta(x) = 0$, unless $\Delta = \emptyset$. Thus, if $T$ is unimodular then the sum in Theorem~\ref{thm:betkemcmullen} collapses to $h_{ \link(\emptyset) } (x) \, B_\emptyset (x) = h_T(x)$, and so Theorem~\ref{thm:unimodularehrharth} is a corollary to Theorem~\ref{thm:betkemcmullen}.
(This argument also shows, in general, that $\h_\P(x) \ge h_T(x)$ component-wise.)
Furthermore, since all ingredients for the sum in Theorem~\ref{thm:betkemcmullen} are nonnegative, this gives another (and the first combinatorial) proof of Theorem~\ref{thm:nonneg}.

Theorem~\ref{thm:betkemcmullen} was greatly extended by Stanley in (and served as some motivation to) his work on local $h$-vectors of subdivisions
\cite{stanleylocalhvectors}; see Athanasiadis' contribution to this volume.
Payne gave a different, multivariate generalization of Theorem~\ref{thm:betkemcmullen} in~\cite{payneehrharttriang}.

Theorem~\ref{thm:betkemcmullen} has a powerful consequence when $\P$ has an interior lattice point; this consequence was fully realized only by Stapledon \cite{stapledondelta} who extended it to general lattice polytopes---Theorem~\ref{thm:stapledonab} below.
Namely, if a lattice polytope $\P$ has an interior lattice point, it admits a regular triangulation that is a cone (at this point) over a boundary triangulation. This has the charming effect that each $h_{ \link(\Delta) } (x)$ appearing in Theorem~\ref{thm:betkemcmullen} is palindromic (due to the afore-mentioned Dehn--Sommerville equations). Since the box polynomials are palindromic and both kinds of polynomials have nonnegative coefficients, a little massaging of the identity in Theorem~\ref{thm:betkemcmullen} gives:

\begin{corollary}\label{cor:betkemcmullen}
Suppose $\P$ is a $d$-dimensional lattice polytope that contains an interior lattice point.
Then there exist unique polynomials $a(x)$ and $b(x)$ with nonnegative coefficients such that
\[
  \h_\P(x) = a(x) + x \, b(x) \, ,
\]
$a(x) = x^d \, a(\frac 1 x)$, and $b(x) = x^{ d-1 } \, b(\frac 1 x)$.
\end{corollary}

The identities for $a(x)$ and $b(x)$ say that $a(x)$ and $b(x)$ are palindromic polynomials; the degree of $a(x)$ is necessarily $d$, while the degree of $b(x)$ is
$d-1$ or smaller; in fact, $b(x)$ can be zero---this happens if and only if $\P$ is the translate of a \emph{reflexive} polytope (i.e., a lattice
polytope whose dual is also a lattice polytope), due to Theorem~\ref{thm:hibi} below. 
Stapledon recently introduced a weighted variant of the Ehrhart $\h$-vector which is
\emph{always} palindromic, motivated by motivic integration and the cohomology of certain toric
varieties \cite{stapledonweightedehrart}. One can easily recover $\h_\P(x)$ from this weighted
Ehrhart $\h$-vector, but one can also deduce the palindromy of both $a(x)$ and $b(x)$ as 
coming from the same source (and this perspective has some serious geometric applications).

The statements that $a(x)$ and $b(x)$ in Corollary~\ref{cor:betkemcmullen} have nonnegative coefficients are straightforward translations of
Hibi's and Stanley's inequalities on Ehrhart $h$-vectors mentioned above (right after and in Theorem~\ref{thm:stanlineq}), in the case that the dimension and the
degree of $\P$ are equal (which is equivalent to $\P$ containing an interior lattice point).
The full generality of Theorem~\ref{thm:stanlineq} as well as Theorem~\ref{thm:stapledonineq} follow from the following generalization of Corollary~\ref{cor:betkemcmullen}.

\begin{theorem}[Stapledon \cite{stapledondelta}]\label{thm:stapledonab}
Suppose $\P$ is a $d$-dimensional lattice polytope of degree $s$ and codegree $l = d+1-s$.
Then there exist unique polynomials $a(x)$ and $b(x)$ with nonnegative coefficients such that
\[
  \left( 1 + x + \dots + x^{ l-1 } \right) \h_\P(x) = a(x) + x^l \, b(x) \, ,
\]
$a(x) = x^d \, a(\frac 1 x)$, $b(x) = x^{ d-l } \, b(\frac 1 x)$, and, writing $a(x) = a_d x^d + a_{ d-1 } x^{ d-1 } + \dots + a_0$,
\[
  1 = a_0 \le a_1 \le a_j \quad \text{ for } \quad 2 \le j \le d-1.
\]
\end{theorem}

Stapledon has recently improved this theorem further, giving infinitely many classes of linear
inequalities among Ehrhart-$h$ coefficients \cite{stapledonadditive}. This exciting new line of
research involves additional techniques from additive number theory. 
 
One can, on the other hand, ask which classes of polytopes satisfy more special sets of equalities or inequalities. Arguably the most natural such equalities/inequalities are those expressing palindromy and unimodality. 

Lattice polytopes with palindromic Ehrhart $h$-vectors are completely classified by the following theorem, which first explicitly surfaced in Hibi's work
on reflexive polytopes but can be traced back to Stanley's work on Hilbert functions of Gorenstein rings.

\begin{theorem}[Hibi--Stanley \cite{hibidual,stanleyhilbertgradedalgebras}]\label{thm:hibi}
If $\P$ is a lattice polytope of degree $s$ and codegree $l = d+1-s$, then its Ehrhart $h$-vector is palindromic if and only if $l \P$ is a translate of a
reflexive polytope.
\end{theorem}

The following result is a start towards a unimodality classification; it was proved by Athanasiadis \cite{athanasiadishstareulerian} and independently by Hibi and Stanley (unpublished).

\begin{theorem}[Athanasiadis--Hibi--Stanley \cite{athanasiadishstareulerian}]\label{thm:athan}
If the $d$-dimensional lattice polytope $\P$ admits a regular unimodular triangulation, then
\[
  \h_{ \lfloor \frac{ d+1 }{ 2 } \rfloor } \ge \dots \ge \h_{ d-1 } \ge \h_d
\]
and
\[
  \h_j \le \binom{ \h_1 + j - 1 }{ j } \quad \text{ for } \quad 0 \le j \le d.
\]
\end{theorem}

Stapledon's work in \cite{stapledondelta} implies further that if the \emph{boundary} of $\P$ admits a regular unimodular triangulation, then
\[
  \h_{ j+1 } \ge \h_{ d-j } \quad \text{ for } \quad 0 \le j \le \lfloor \tfrac d 2 \rfloor - 1 
\]
(which was also proved in \cite{athanasiadishstareulerian} under the stronger assumption that $\P$ admits a regular unimodular triangulation) and
\[
  \h_0 + \dots + \h_{ j+1 } \le \h_d + \dots + \h_{ d-j } + \binom{ \h_1 - \h_d + j + 1 }{ j+1 } \quad \text{ for } \quad 0 \le j \le \lfloor \tfrac d 2 \rfloor - 1.
\]
Naturally, if in addition to the conditions in Theorem~\ref{thm:athan}, $\P$ has degree $d$ and $\h_\P(x)$ is palindromic, then $\h_\P(x)$ is unimodal.
The proof of Theorem~\ref{thm:athan} starts with Theorem~\ref{thm:unimodularehrharth} and then shows that the $h$-vector of the unimodular triangulation satisfies
the stated inequalities.
Athanasiadis' approach was inspired by work of Reiner and Welker on order polytopes of graded posets and a connection between the \emph{Charney--Davis} and
\emph{Neggers--Stanley conjectures} \cite{reinerwelker} and can be taken further: Athanasiadis used similar methods to prove Stanley's conjecture mentioned in
Section~\ref{sec:anand} that the Ehrhart $h$-vector of the Birkhoff--von Neumann polytope is unimodal \cite{athanasiadismagic}. Bruns and R\"omer generalized this to
any Gorenstein polytope that admits a regular unimodular triangulation~\cite{brunsroemer}.

Going into a somewhat different direction, Schepers and van Langenhoven recently proved that lattice parallelepipeds have a unimodal Ehrhart $h$-vector
\cite{schepersvanlangenhoven}. The conjecture that any integrally closed polytope (of which both parallelepipeds and the Birkhoff--von Neumann polytope are examples) has a unimodal Ehrhart $h$-vector remains open, even for integrally closed reflexive polytopes (though recent work of Braun and Davis give some pointers of what could be tried here~\cite{braundavis}); this is closely related to a conjecture of Stanley that every standard graded Cohen--Macaulay domain has a unimodal $h$-vector~\cite{stanleylogconcave}.


\bibliographystyle{amsplain}
\bibliography{bib}

\def\cprime{$'$} \def\cprime{$'$}
\providecommand{\bysame}{\leavevmode\hbox to3em{\hrulefill}\thinspace}
\providecommand{\MR}{\relax\ifhmode\unskip\space\fi MR }
\providecommand{\MRhref}[2]{%
  \href{http://www.ams.org/mathscinet-getitem?mr=#1}{#2}
}
\providecommand{\href}[2]{#2}
\begin{thebibliography}{10}

\bibitem{ananddumirgupta}
Harsh Anand, Vishwa~Chander Dumir, and Hansraj Gupta, \emph{A combinatorial
  distribution problem}, Duke Math. J. \textbf{33} (1966), 757--769.

\bibitem{andrewspa1}
George~E. Andrews, \emph{Mac{M}ahon's partition analysis. {I}. {T}he lecture
  hall partition theorem}, Mathematical essays in honor of {G}ian-{C}arlo
  {R}ota ({C}ambridge, {MA}, 1996), Progr. Math., vol. 161, Birkh\"auser
  Boston, Boston, MA, 1998, pp.~1--22.

\bibitem{athanasiadishstareulerian}
Christos~A. Athanasiadis, \emph{{$h^\ast$}-vectors, {E}ulerian polynomials and
  stable polytopes of graphs}, Electron. J. Combin. \textbf{11} (2004/06),
  no.~2, Research Paper 6, 13 pp. (electronic).

\bibitem{athanasiadismagic}
\bysame, \emph{Ehrhart polynomials, simplicial polytopes, magic squares and a
  conjecture of {S}tanley}, J. Reine Angew. Math. \textbf{583} (2005),
  163--174, {\tt arXiv:math/0312031}.

\bibitem{batyrevhvector}
Victor~V. Batyrev, \emph{Lattice polytopes with a given {$h\sp *$}-polynomial},
  Algebraic and geometric combinatorics, Contemp. Math., vol. 423, Amer. Math.
  Soc., Providence, RI, 2006, pp.~1--10, {\tt arXiv:math.CO/0602593}.

\bibitem{beckpixton}
Matthias Beck and Dennis Pixton, \emph{The {E}hrhart polynomial of the
  {B}irkhoff polytope}, Discrete Comput. Geom. \textbf{30} (2003), no.~4,
  623--637, {\tt arXiv:math.CO/0202267}.

\bibitem{ccd}
Matthias Beck and Sinai Robins, \emph{Computing the {C}ontinuous {D}iscretely:
  Integer-point {E}numeration in {P}olyhedra}, Undergraduate Texts in
  Mathematics, Springer, New York, 2007, electronically available at {\tt
  http://math.sfsu.edu/beck/ccd.html}.

\bibitem{bsstanleyirrational}
Matthias Beck and Frank Sottile, \emph{Irrational proofs for three theorems of
  {S}tanley}, European J. Combin. \textbf{28} (2007), no.~1, 403--409, {\tt
  arXiv:math.CO/0506315}.

\bibitem{nnz}
Matthias Beck and Thomas Zaslavsky, \emph{The number of nowhere-zero flows on
  graphs and signed graphs}, J. Combin. Theory Ser. B \textbf{96} (2006),
  no.~6, 901--918, {\tt arXiv:math.CO/0309331}.

\bibitem{betkemcmullen}
Ulrich Betke and Peter McMullen, \emph{Lattice points in lattice polytopes},
  Monatsh. Math. \textbf{99} (1985), no.~4, 253--265.

\bibitem{braundavis}
Benjamin Braun and Robert Davis, \emph{Ehrhart series, unimodality, and
  integrally closed reflexive polytopes}, Preprint ({\tt
  arXiv:math/1403.5378}), 2014.

\bibitem{brunssemimagic}
Winfried Bruns, \emph{Commutative algebra arising from the
  {A}nand--{D}umir--{G}upta conjectures}, Commutative algebra and
  combinatorics, Ramanujan Math. Soc. Lect. Notes Ser., vol.~4, Ramanujan Math.
  Soc., Mysore, 2007, pp.~1--38.

\bibitem{brunsgubeladzektheory}
Winfried Bruns and Joseph Gubeladze, \emph{Polytopes, rings, and {$K$}-theory},
  Springer Monographs in Mathematics, Springer, Dordrecht, 2009.

\bibitem{brunsroemer}
Winfried Bruns and Tim R{\"o}mer, \emph{{$h$}-vectors of {G}orenstein
  polytopes}, J. Combin. Theory Ser. A \textbf{114} (2007), no.~1, 65--76.

\bibitem{canfieldmckay}
E.~Rodney Canfield and Brendan~D. McKay, \emph{The asymptotic volume of the
  {B}irkhoff polytope}, Online J. Anal. Comb. (2009), no.~4, 4 pages, {\tt
  arXiv:math.0705.2422}.

\bibitem{deloeraliuyoshida}
Jes{\'u}s~A. De~Loera, Fu~Liu, and Ruriko Yoshida, \emph{A generating function
  for all semi-magic squares and the volume of the {B}irkhoff polytope}, J.
  Algebraic Combin. \textbf{30} (2009), no.~1, 113--139.

\bibitem{deloerarambausantos}
Jes{\'u}s~A. De~Loera, J{\"o}rg Rambau, and Francisco Santos,
  \emph{Triangulations}, Algorithms and Computation in Mathematics, vol.~25,
  Springer-Verlag, Berlin, 2010.

\bibitem{ehrhartpolynomial}
Eug{\`e}ne Ehrhart, \emph{Sur les poly\`edres rationnels homoth\'etiques \`a
  {$n$}\ dimensions}, C. R. Acad. Sci. Paris \textbf{254} (1962), 616--618.

\bibitem{eisenbudbook}
David Eisenbud, \emph{Commutative {A}lgebra {W}ith a {V}iew {T}oward
  {A}lgebraic {G}eometry}, Graduate Texts in Mathematics, vol. 150,
  Springer-Verlag, New York, 1995.

\bibitem{elliott}
E.~B. Elliott, \emph{On linear homogeneous diophantine equations}, Quartely J.
  Pure Appl. Math. \textbf{34} (1903), 348--377.

\bibitem{haasenillpayne}
Christian Haase, Benjamin Nill, and Sam Payne, \emph{Cayley decompositions of
  lattice polytopes and upper bounds for {$h^*$}-polynomials}, J. Reine Angew.
  Math. \textbf{637} (2009), 207--216, {\tt arXiv:math/0804.3667}.

\bibitem{hibiehrhartineq}
Takayuki Hibi, \emph{Some results on {E}hrhart polynomials of convex
  polytopes}, Discrete Math. \textbf{83} (1990), no.~1, 119--121.

\bibitem{hibidual}
\bysame, \emph{Dual polytopes of rational convex polytopes}, Combinatorica
  \textbf{12} (1992), no.~2, 237--240.

\bibitem{hibilowerbound}
\bysame, \emph{A lower bound theorem for {E}hrhart polynomials of convex
  polytopes}, Adv. Math. \textbf{105} (1994), no.~2, 162--165.

\bibitem{hochster}
Melvin Hochster, \emph{Rings of invariants of tori, {C}ohen--{M}acaulay rings
  generated by monomials, and polytopes}, Ann. of Math. (2) \textbf{96} (1972),
  318--337.

\bibitem{lagariasziegler}
Jeffrey~C. Lagarias and G{\"u}nter~M. Ziegler, \emph{Bounds for lattice
  polytopes containing a fixed number of interior points in a sublattice},
  Canad. J. Math. \textbf{43} (1991), no.~5, 1022--1035.

\bibitem{macdonald}
Ian~G. Macdonald, \emph{Polynomials associated with finite cell-complexes}, J.
  London Math. Soc. (2) \textbf{4} (1971), 181--192.

\bibitem{macmahon}
Percy~A. MacMahon, \emph{Combinatory {A}nalysis}, Chelsea Publishing Co., New
  York, 1960, reprint of the 1915 original.

\bibitem{payneehrharttriang}
Sam Payne, \emph{Ehrhart series and lattice triangulations}, Discrete Comput.
  Geom. \textbf{40} (2008), no.~3, 365--376, {\tt arXiv:math/0702052}.

\bibitem{pick}
Georg~Alexander Pick, \emph{Geometrisches zur {Z}ahlenlehre}, Sitzenber.~Lotos
  (Prague) \textbf{19} (1899), 311--319.

\bibitem{reinerwelker}
Victor Reiner and Volkmar Welker, \emph{On the {C}harney--{D}avis and
  {N}eggers--{S}tanley conjectures}, J. Combin. Theory Ser. A \textbf{109}
  (2005), no.~2, 247--280.

\bibitem{schepersvanlangenhoven}
Jan Schepers and Leen Van~Langenhoven, \emph{Unimodality questions for
  integrally closed lattice polytopes}, Ann. Comb. \textbf{17} (2013), no.~3,
  571--589, {\tt arXiv:math/1110.3724}.

\bibitem{scott}
Paul~R. Scott, \emph{On convex lattice polygons}, Bull. Austral. Math. Soc.
  \textbf{15} (1976), no.~3, 395--399.

\bibitem{stanleythesis}
Richard~P. Stanley, \emph{Ordered structures and partitions}, American
  Mathematical Society, Providence, R.I., 1972, Memoirs of the American
  Mathematical Society, No. 119.

\bibitem{stanleymagic}
\bysame, \emph{Linear homogeneous {D}iophantine equations and magic labelings
  of graphs}, Duke Math. J. \textbf{40} (1973), 607--632.

\bibitem{stanleyreciprocity}
\bysame, \emph{Combinatorial reciprocity theorems}, Advances in Math.
  \textbf{14} (1974), 194--253.

\bibitem{stanleyupperbound}
\bysame, \emph{The upper bound conjecture and {C}ohen--{M}acaulay rings},
  Studies in Appl. Math. \textbf{54} (1975), no.~2, 135--142.

\bibitem{stanleymagiccohenmac}
\bysame, \emph{Magic labelings of graphs, symmetric magic squares, systems of
  parameters, and {C}ohen--{M}acaulay rings}, Duke Math. J. \textbf{43} (1976),
  no.~3, 511--531.

\bibitem{stanleyhilbertgradedalgebras}
\bysame, \emph{Hilbert functions of graded algebras}, Advances in Math.
  \textbf{28} (1978), no.~1, 57--83.

\bibitem{stanleydecomp}
\bysame, \emph{Decompositions of rational convex polytopes}, Ann. Discrete
  Math. \textbf{6} (1980), 333--342.

\bibitem{stanleylogconcave}
\bysame, \emph{Log-concave and unimodal sequences in algebra, combinatorics,
  and geometry}, Graph theory and its applications: {E}ast and {W}est ({J}inan,
  1986), Ann. New York Acad. Sci., vol. 576, New York Acad. Sci., New York,
  1989, pp.~500--535.

\bibitem{stanleyinequ}
\bysame, \emph{On the {H}ilbert function of a graded {C}ohen--{M}acaulay
  domain}, J. Pure Appl. Algebra \textbf{73} (1991), no.~3, 307--314.

\bibitem{stanleylocalhvectors}
\bysame, \emph{Subdivisions and local {$h$}-vectors}, J. Amer. Math. Soc.
  \textbf{5} (1992), no.~4, 805--851.

\bibitem{stanleymonotonicity}
\bysame, \emph{A monotonicity property of {$h$}-vectors and {$h^*$}-vectors},
  European J. Combin. \textbf{14} (1993), no.~3, 251--258.

\bibitem{stanleycombcommalg}
\bysame, \emph{Combinatorics and {C}ommutative {A}lgebra}, second ed., Progress
  in Mathematics, vol.~41, Birkh\"auser Boston Inc., Boston, MA, 1996.

\bibitem{stanleyec1}
\bysame, \emph{Enumerative {C}ombinatorics. {V}olume 1}, second ed., Cambridge
  Studies in Advanced Mathematics, vol.~49, Cambridge University Press,
  Cambridge, 2012.

\bibitem{stanleyhowupperbound}
\bysame, \emph{How the upper bound conjecture was proved}, {\tt
  http://www-math.mit.edu/$\sim$rstan/papers.html}, 2013.

\bibitem{stapledonweightedehrart}
Alan Stapledon, \emph{Weighted {E}hrhart theory and orbifold cohomology}, Adv.
  Math. \textbf{219} (2008), no.~1, 63--88, {\tt arXiv:math/0711.4382}.

\bibitem{stapledondelta}
\bysame, \emph{Inequalities and {E}hrhart {$\delta$}-vectors}, Trans. Amer.
  Math. Soc. \textbf{361} (2009), no.~10, 5615--5626, {\tt
  arXiv:math/0801.0873}.

\bibitem{stapledonadditive}
\bysame, \emph{Additive number theory and inequalities in {E}hrhart theory},
  Preprint ({\tt arXiv:0904.3035v2}), 2010.

\bibitem{xinmonsterrec}
Guoce Xin, \emph{Generalization of {S}tanley's monster reciprocity theorem}, J.
  Combin. Theory Ser. A \textbf{114} (2007), no.~8, 1526--1544, {\tt
  arXiv:math/0504425}.

\bibitem{zaslavskyorientationsignedgraphs}
Thomas Zaslavsky, \emph{Orientation of signed graphs}, European J. Combin.
  \textbf{12} (1991), no.~4, 361--375.

\bibitem{ziegler}
G{\"u}nter~M. Ziegler, \emph{Lectures on {P}olytopes}, Springer-Verlag, New
  York, 1995.

\end{thebibliography}

\end{document}